\documentclass{amsart}
\usepackage{amsfonts}
\usepackage{amsmath}
\usepackage{enumerate}
\usepackage{amsthm}
\usepackage{hyperref}
\usepackage{cite}

\begin{document}

\newtheorem{definition}{Definition}[section]
\newtheorem{theorem}[definition]{Theorem}
\newtheorem{proposition}[definition]{Proposition}
\newtheorem{remark}[definition]{Remark}
\newtheorem{lemma}[definition]{Lemma}
\newtheorem{corollary}[definition]{Corollary}
\newtheorem{example}[definition]{Example}

\numberwithin{equation}{section}

\title[The Kazdan-Warner problem...]{The Kazdan-Warner problem on compact K\"ahler surfaces}
\author{Weike Yu}

\date{}

\begin{abstract}
In this paper, we investigate a Kazdan-Warner problem on compact K\"ahler surfaces, which corresponds to prescribing sign-changing Chern scalar curvatures, and establish a Chen-Li type existence theorem on compact K\"ahler surfaces when the candidate curvature function is of negative average. Moreover, we give an alternative proof of Ding-Liu's theorem [Trans. Amer. Math. Soc. 347(1995) 1059-1066] on prescribing sign-changing Gaussian curvatures.
\end{abstract}
\maketitle
\section{Introduction}
Suppose that $M$ is a two-dimensional compact Riemannian manifold with a Riemannian metric $g$. A classical problem raised by Kazdan and Warner is the following question: given a smooth function $K: M\rightarrow \mathbb{R}$, can it be realized as the Gaussian curvature of some metric $\hat{g}$ which is conformal to $g$? Let $\hat{g}=e^{2u}g$, then this problem is equivalent to the solvability of the semilinear elliptic equation
\begin{align}\label{1.1}
-\Delta_g u+k=Ke^{2u},
\end{align}
where $\Delta_g$ is the Laplacian with respect to $g$ and $k$ is the Gaussian curvature of $(M, g)$. 

According to the Uniformization Theorem, without loss of generality, one can assume that $(M, g)$ has constant Gaussian curvature $k\in \mathbb{R}$.
In \cite{[KW]}, Kazdan and Warner studied the above equation on $(M, g)$ with number $k<0$, which corresponds to the prescribing Gaussian curvature problem on compact Riemannian surfaces with negative Euler characteristic. They showed that
\begin{enumerate}[(i)]
\item If the equation \eqref{1.1} has a smooth solution, then $\int_M K d\mu_g<0$.
\item If $\int_M K d\mu_g<0$, then there exists a constant $k_*\in [-\infty,0)$ such that \eqref{1.1} has a smooth solution if $0>k>k_*$, but does not have a smooth solution if $k<k_*$.
\item $k_*=-\infty$ if and only if $K\leq 0\ (\not\equiv0)$.
\end{enumerate}
Furthermore, for the case $k=k_*>-\infty$, W. Chen and C. Li \cite{[CL]} proved that
\begin{enumerate}[(iv)]
\item Let $k_*>-\infty$ be defined above, then \eqref{1.1} has at least one smooth solution when $k=k_*$.
\end{enumerate}
Note that in \cite[Sect.10]{[KW]}, Kazdan-Warner actually proved that the above conclusions (i)-(iii) hold for Eq. \eqref{1.1} with number $k<0$ on any compact Riemannian manifolds with arbitrary dimension. Besides, from (i) and (iii), it follows that $k_*>-\infty$ if and only if $K$ is sign-changing and $\int_M K d\mu_g<0$.

In this paper, we aim to generalize the result of W. Chen and C. Li \cite{[CL]} to higher-dimensional case. Let $(M^n, \omega)$ be a compact Hermitian manifold, where $\omega$ is a Hermitian metric and $\dim_{\mathbb{C}}M=n\geq 2$. The prescribed Chern scalar curvature problem (see Section 2 for details) is a Hermitian analogue of prescribing scalar curvatures on Riemannian manifolds: given a smooth function $S: M\rightarrow\mathbb{R}$, does there exist a Hermitian metric $\hat{\omega}$ conformal to $\omega$, that is, $\hat{\omega}=e^{\frac{2}{n}u}\omega$, such that its Chern scalar curvature is $S$? This problem is equivalent to the solvability of 
\begin{align}\label{1.2.}
-\Delta^{Ch}_\omega u+S^{Ch}(\omega) =S e^{\frac{2}{n}u},
\end{align}
where $\Delta^{Ch}_\omega$ is the Chern Laplacian with respect to $\omega$ and $S^{Ch}(\omega)$ is the Chern scalar curvature of $(M, \omega)$. When $S$ is a constant, the above problem is referred to as the Chern-Yamabe problem, which was first proposed by Angella-Calamai-Spotti in \cite{[ACS]}.

If the Gauduchon degree of $\{\omega\}$ (see \eqref{2.2}) is negative, according to \cite[Theorem 4.1]{[ACS]}, it is natural to consider the case that $S^{Ch}(\omega)$ in \eqref{1.2.} is a negative constant. Now we free equation \eqref{1.2.} from the geometric situation and instead consider the following equation:
\begin{align}\label{1.2}
-\Delta^{Ch}_\omega u+\alpha =S e^{\frac{2}{n}u},
\end{align}
on arbitrary compact Hermitian manifold $(M^n, \omega)$ without any assumptions on curvatures, where $\alpha$ is a negative constant and $S\in C^\infty(M)$. In \cite{[Fus]}, Fusi proved that
\begin{enumerate}[(1)]
\item If the equation \eqref{1.2} has a smooth solution, then $\int_M Sf_0d\mu_\omega<0$.
\item If $\int_M Sf_0d\mu_\omega<0$, then there exists a constant $\alpha_*\in [-\infty,0)$ depending on $S$ such that one can solve \eqref{1.2} for all $\alpha\in (\alpha_*,0)$, but cannot solve \eqref{1.2} if $\alpha<\alpha_*$.
\item $\alpha_*=-\infty$ if $S\leq 0\ (\not\equiv0)$.
\end{enumerate}
Here the function $f_0\in C^\infty(M)$ is the eccentricity function associated to $\omega$, namely, $f_0\in \text{ker}((\Delta^{Ch}_\omega)^*)$ and $\langle f_0,1\rangle_{L^2(M)}=\text{Vol(M)}$, where $(\Delta^{Ch}_\omega)^*$ is the formal adjoint of $\Delta^{Ch}_\omega$. In \cite{[Gau1]}, Gauduchon proved that $f_0>0$, and $f_0=1$ if and only if $\omega$ is Gauduchon, i.e., $\partial \bar{\partial} w^{n-1}=0$. Furthermore, following the same argument as in \cite[Theorem 10.5]{[KW]}, one can also obtain
\begin{enumerate}[(3')]
\item $\alpha_*=-\infty$ if and only if $S\leq 0\ (\not\equiv0)$.
\end{enumerate}

In this paper, we consider the case $\alpha=\alpha_*>-\infty$, by the spirit of our recent paper \cite{[Yu1]}, we establish the following result:
\begin{theorem}\label{thm1.1}
Let $(M, \omega)$ be a compact K\"ahler surface and let $S$ be a smooth function on $M$ with $\int_M S d\mu_\omega<0$. Suppose that $\alpha_*>-\infty$, then \eqref{1.2} has at least one solution for $\alpha=\alpha_*$.
\end{theorem}
\begin{remark}
As mentioned above, Chen-Li \cite{[CL]} proved the above result on the Riemannian surface (see (iv) in the present section). Their proof relies on the estimates in Lemma \ref{lem4.1}, which, however, seems difficult to generalize to higher-dimensional cases, so their method cannot be used to prove the above theorem. In addition, we also give a different proof of Chen-Li's result in Section $3$. 
\end{remark}
\begin{remark}
At the end of this paper (see Appendix), we give a simple proof of the key local $C^0$-estimate in the proof of Theorem \ref{thm1.1} by using the classical maximum principle.
\end{remark}

Finally, in terms of the $\sup+\inf$ inequality \cite{[BLS]}, we give an alternative proof of the following Ding-Liu's theorem \cite{[DL]} for prescribing sign-changing Gaussian curvatures on compact Riemannian surfaces.

\begin{theorem}[cf. \cite{[DL]}]
Let $(M,\omega)$ be a compact Riemannian surfaces with Euler characteristic $\chi (M)<0$. Let $g_0$ be a nonconstant smooth function on $M$ with $\max_M g_0=0$. Then there exists a constant $\lambda^*\in (0,-\min_M g_0)$ such that 
\begin{enumerate}
\item If $\lambda\in (0, \lambda^*)$, then $g_0+\lambda\in PC(\omega)$;\\
 If $\lambda\in (\lambda^*,+\infty)$, then $g_0+\lambda\not\in PC(\omega)$.
 \item If $\lambda=\lambda^*$,  then $g_0+\lambda^*\in PC(\omega)$.
\end{enumerate}
Here $PC(\omega)$ denotes the set of $C^\infty(M)$ functions which are Gaussian curvatures of all $\tilde\omega\in \{\omega\}$.
\end{theorem}

This paper is organized as follows. In Section 2, we recall some basic notions and notations related to the prescribed Chern scalar curvature problem. In Section 3, we prove our main result (Theorem \ref{thm1.1}). In Section 4, we give an alternative proof of Ding-Liu's theorem on prescribing sign-changing Gaussian curvatures. In Appendix, we give a simple proof of the key local $C^0$-estimate in the proof of Theorem \ref{thm1.1}.

\textbf{Acknowledgements.} The author would like to thank Prof. Yuxin Dong and Prof. Xi Zhang for their continued support and encouragement.

\section{The prescribed Chern scalar curvature problem}
In this section, we will give a brief introduction to the prescribed Chern scalar curvature problem on compact Hermitian manifolds (cf. \cite{[ACS], [CZ], [Fus], [Ho], [LM], [Yu1], [Yu2]}, etc.). 

Let $(M^n, J, h)$ be a Hermitian manifold with complex dimension $n$, where $h$ is a Hermitian metric and $J$ is a complex structure. The corresponding fundamental form $\omega$ of $h$ is given by $\omega(\cdot,\cdot)=h(J\cdot, \cdot)$, and we will confuse the Hermitian metric $h$ and its fundamental form $\omega$ in this paper. It is well-known that there exists a unique affine connection $\nabla^{Ch}$ preserving both the Hermitian metric $h$ and the complex structure $J$, that is, $\nabla^{Ch}h=0$, $\nabla^{Ch}J=0$, whose torsion $T^{Ch}(X,Y)$ satisfies $T^{Ch}(JX,Y)=T^{Ch}(X,JY)$ for any tangent vectors $X, Y\in TM$. Such a connection $\nabla^{Ch}$ is called the Chern connection.

For a Hermitian manifold $(M^n, \omega)$ with Chern connection $\nabla^{Ch}$, the corresponding Chern scalar curvature is given by
\begin{align}
S^{Ch}(\omega)=\text{tr}_{\omega}Ric^{(1)}(\omega)=\text{tr}_{\omega}\sqrt{-1}\bar{\partial}\partial \log{\omega^n},
\end{align}
where $Ric^{(1)}(\omega)$ is the first Chern Ricci curvature. 

Let $\{\omega\}=\{e^{\frac{2}{n}u}\omega\ |\ u\in C^\infty(M)\}$ denote the conformal class of the Hermitian metric $\omega$. In \cite{[Gau1]}, Gauduchon proved that
\begin{theorem}\label{theorem2.2}
Let $(M^n, \omega)$ be a compact Hermitian manifold with complex dimension $n\geq 2$. Then there exists a unique Gauduchon metric $\eta\in \{\omega\}$ (i.e., $\partial \bar{\partial} w^{n-1}=0$) with volume $1$.
\end{theorem}
Using the above theorem, one can define an invariant $\Gamma(\{\omega\})$ of the conformal class $\{\omega\}$ as follows:
\begin{align}\label{2.2}
\Gamma(\{\omega\})=\frac{1}{(n-1)!}\int_M c^{BC}_1(K^{-1}_M)\wedge\eta^{n-1}=\int_M S^{Ch}(\eta)d\mu_\eta,
\end{align}
where $c^{BC}_1(K^{-1}_M)$ is the first Bott-Chern class of anti-canonical line bundle $K^{-1}_M$, and $d\mu_\eta$ denotes the volume form of the Gauduchon metric $\eta$.

Consider the conformal change $\widetilde{\omega}=e^{\frac{2}{n}u}\omega$ on the compact Hermitian manifold $(M^n, \omega)$. From \cite{[Gau]}, the Chern scalar curvatures of $\widetilde{\omega}$ and $\omega$ have the following relationship:
\begin{align}\label{2.7}
-\Delta^{Ch}_\omega u+S^{Ch}(\omega)=S^{Ch}(\widetilde{\omega})e^{\frac{2}{n}u},
\end{align}
where $\Delta^{Ch}_\omega$ is the Chern Laplacian defined by $\Delta^{Ch}u=-2\sqrt{-1} tr_{\omega}\overline{\partial}\partial u$ for any $u\in C^\infty(M)$. Note that $-\Delta^{Ch}_\omega u=-\Delta_d u+(du,\theta)_{\omega}$, where $\Delta_du=-d^*du$ is the Hodge-de Rham Laplacian, $\theta$ is the Lee form given by $d\omega^{n-1}=\theta\wedge \omega^{n-1}$ and $(\cdot, \cdot)_\omega$ denotes the inner product on $1$-form induced by $\omega$. Hence, If $\omega$ is K\"ahler or balanced, then $\Delta^{Ch}_\omega u=\Delta_d u$ for any smooth function $u$.

Finally, we recall two properties about classical Sobolev spaces on the compact Hermitian manifold $(M^n,\omega)$, which will be used later. 
\begin{lemma}\label{lem2.2}
Assume that $n\geq2$. Then there exists a uniform constant $C$ such that 
\begin{align}
\left(\int_M|f|^{2\beta}d\mu_\omega \right)^{\frac{1}{\beta}}\leq C\left( \int_M|\nabla f|^{2}d\mu_\omega+\int_M|f|^{2}d\mu_\omega\right)
\end{align}
for any $f\in W^{1,2}(M)$, where $\beta=\frac{n}{n-1}$ and $W^{1,2}(M)$ is the classical Sobolev space on $(M^n, \omega)$.
\end{lemma}

\begin{lemma}[\cite{[Fon]}]\label{lemma2.4}
For any $f\in W^{k,p}(M)$ with $\int_M fd\mu_\omega=0$, $\int_M |\nabla^k f|^{p}d\mu_\omega\leq1$, $kp=2n$, $k\in \mathbb{N+}$, $p>1$, there exist two constants $\delta_1=\delta_1(k,n)$ and $\delta_2=\delta_2(k, M)$ such that
\begin{align}
\int_M e^{\delta_1|f|^{\frac{p}{p-1}}} d\mu_\omega\leq \delta_2.
\end{align}
\end{lemma}

\section{Proof of Theorem \ref{thm1.1}}

In this section, we will give a proof of Theorem \ref{thm1.1}. Consider the following equation on the compact K\"ahler manifold $(M^n, \omega)$ for $n\geq 1$:
\begin{align}\label{3.1}
-\Delta^{Ch}_\omega u+\alpha_k=S e^{\frac{2}{n}u},
\end{align}
where $\{\alpha_k\}$ is a sequence of numbers with $\alpha_*<\alpha_k<\alpha^*<0$ and $\lim_{k\rightarrow +\infty}\alpha_k=\alpha_*$, $\alpha^*<0$ is a fixed constant. Since $\alpha_k>\alpha_*$, there always exists a solution $u_k\in C^\infty(M)$ of \eqref{3.1}. In general, the solution $u_k$ is not unique.

For each $\alpha_k$, we choose $\tilde{\alpha}_k$ satisfying $\alpha_*<\tilde{\alpha}_k<\alpha_k$. By the fact (2) in Section 1, there exists a solution $u_+\in C^\infty(M)$ of 
\begin{align}
-\Delta^{Ch}_\omega u_+ +\tilde{\alpha}_k=S e^{\frac{2}{n}u_+},
\end{align}
which implies that $u_+$ is a super-solution of \eqref{3.1}, that is
\begin{align}
-\Delta^{Ch}_\omega u_+ +\alpha_k>S e^{\frac{2}{n}u_+}.
\end{align}
Let $u_-=-C_k<0$ be a constant with $C_k>\max\{-\frac{n}{2}\log (\frac{\alpha^*}{\inf_M S}), -\inf_M u_+\}$, then it satisfies
\begin{align}
-\Delta^{Ch}_\omega u_-+\alpha_k<S e^{\frac{2}{n}u_-}\quad \text{and}\quad u_-<u_+.
\end{align}
Here both $u_-$ and $u_+$ are related to $k$. Through some standard variational argument (cf. \cite[Lemma 4.3]{[Yu1]}), we can conclude that

\begin{lemma}\label{lemma3.1}
There exists a function $u_k\in C^\infty(M)\cap X$ such that 
\begin{align}\label{3.5...}
I_k(u_k)=\inf_{u\in X}I_k(u)=\inf_{u\in X}\int_M \left(|\nabla u|^2+2\alpha_ku-nSe^{\frac{2}{n}u}\right)d\mu_\omega,
\end{align}
where $X=\{u\in W^{1,2}(M): u_-\leq u\leq u_+\ \text{a.e.\ on}\ M\}$. (Note that $u_-<u_k<u_+$ because of the maximum principle.) Therefore,
\begin{align}
&0=\frac{d}{dt}|_{t=0}I_k(u_k+t\varphi)=2\int_M\left(\nabla u_k\cdot \nabla\varphi+\alpha_k\varphi-Se^{\frac{2}{n}u_k}\varphi\right)d\mu_\omega,\label{3.3}\\
&0\leq\frac{d^2}{dt^2}|_{t=0}I_k(u_k+t\varphi)=2\int_M|\nabla \varphi|^2d\mu_\omega-\frac{4}{n}\int_M Se^{\frac{2}{n}u_k}\varphi^2 d\mu_\omega\label{3.4}
\end{align}
for any $\varphi\in W^{1,2}(M)$. Hence, by \eqref{3.3}, we obtain that $u_k$ satisfies
\begin{align}\label{3.8...}
-\Delta^{Ch}_\omega u_k+\alpha_k=S e^{\frac{2}{n}u_k}.
\end{align}
\end{lemma}
\begin{remark}
According to \cite[Prop. 5.3]{[ACS]} and \cite[Prop. 2.12]{[Fus]}, we know that \eqref{3.8...} is the Euler-Lagrange equation for the functional defined in \eqref{3.5...} if and only if $\omega$ is balanced (which is K\"ahler in complex dimension $2$). That's the reason why the proof of Theorem \ref{thm1.1} works only for K\"ahler surfaces and not for any Hermitian complex surfaces in general.
\end{remark}

\begin{lemma}\label{lemma 3.3}
There is a constant $A>0$ independent of $k$ such that $u_k>-A$ for any $k=1,2,\dots$.
\end{lemma}
\proof 
Since $\inf_M S<0$ and $\alpha^*<0$, there is a unique solution $f\in C^\infty(M)$ of
\begin{align}
-\Delta^{Ch}_\omega f+\alpha^*=(\inf_M S) e^{\frac{2}{n}f},
\end{align}
which is due to the fact (3) in Section 1. Set $\phi_t=f-t$, where $t\geq 0$, then 
\begin{equation}\label{3.9}
\begin{aligned}
-&\Delta^{Ch}_\omega \phi_t+\alpha_k-S e^{\frac{2}{n}\phi_t}\\
&=-\alpha^*+\alpha_k+(\inf_M S) e^{\frac{2}{n}f}-S e^{\frac{2}{n}f}e^{-\frac{2}{n}t}\\
&<(\inf_M S) e^{-\frac{2}{n}\|f\|_{C^0(M)}}(1-e^{-\frac{2}{n}t})\\
&\leq 0,
\end{aligned}
\end{equation}
for any $t\geq 0$, since $\alpha_k<\alpha^*$ and $\inf_M S<0$. Now we will prove $u_k\geq \phi_0=f$ on $M$. If not, there exists a number $t_0>0$ such that $u_k\geq \phi_{t_0}$ on $M$ and $u_k(x_0)=\phi_{t_0}(x_0)$ for some $x_0\in M$. Then we have at $x_0$
\begin{equation}
\begin{aligned}
-&\Delta^{Ch}_\omega \phi_{t_0}+\alpha_k-S e^{\frac{2}{n}\phi_{t_0}}\\
&=\Delta^{Ch}_\omega (u_k-\phi_{t_0})+(-\Delta^{Ch}_\omega u_k+\alpha_k-S e^{\frac{2}{n}u_k})+S(e^{\frac{2}{n}u_k}-e^{\frac{2}{n}\phi_{t_0}})\\
&=\Delta^{Ch}_\omega (u_k-\phi_{t_0})\geq 0,
\end{aligned}
\end{equation}
because of the maximum principle. This leads to a contradiction with \eqref{3.9}. Therefore, $u_k>-A=\inf_Mf-1$.

\qed

Since $\int_M Sd\mu_\omega<0$, there is a constant $\epsilon_0>0$ such that $M_-=\{x\in M: S(x)<-\epsilon_0\}$ is not empty.

\begin{lemma}\label{lem3.2}
$u_k$ is locally uniformly $W^{1,2}$-bounded in the open set $M_-$ for $n\geq 1$ as $k\rightarrow +\infty$.
\end{lemma}
\proof
Set $v_k=u_k+A>0$, it follows from \eqref{3.3} that
\begin{align}\label{3.5}
\int_{M}\left(\nabla v_k\cdot\nabla \varphi+\alpha_k\varphi-Se^{-\frac{2}{n}A}e^{\frac{2}{n}v_k}\varphi\right)d\mu_\omega=0.
\end{align}
Let $D_1, D_2$ be two arbitrary open subsets of $M_-$ with $D_1\subset\subset D_2\subset\subset M_-$. According to \cite[Lemma 3.3]{[Doa]}, we choose a cut-off function $\xi\in C^\infty(M)$ satisfying
\begin{align}\label{3.12.}
\xi(x)=
\begin{cases}
1& x\in D_1\\
0& x\in M\setminus D_2
\end{cases}
,\ 0\leq\xi\leq1,\ |\nabla \xi|\leq C\xi^{\frac{1}{2}},
\end{align}
where $C$ is a constant. Taking the text function $\varphi=\xi^2v_k$ in \eqref{3.5} yields
\begin{align}
\int_{D_2}\left(\nabla v_k\cdot\nabla (\xi^2v_k)+\alpha_k\xi^2v_k-Se^{-\frac{2}{n}A}e^{\frac{2}{n}v_k}\xi^2v_k\right)d\mu_\omega=0.
\end{align}
Making use of
\begin{align}
\nabla v_k\cdot\nabla (\xi^2v_k)=|\nabla (\xi v_k)|^2-|\nabla \xi|^2v_k^2,\quad S(x)\leq -\epsilon_0<0
\end{align}
for any $x\in D_2$ and \eqref{3.12.}, we deduce that
\begin{equation}
\begin{aligned}
0&\geq\int_{D_2}\left(|\nabla (\xi v_k)|^2-|\nabla \xi|^2v_k^2+\alpha_k\xi^2v_k+\epsilon_0e^{-\frac{2}{n}A}e^{\frac{2}{n}v_k}\xi^2v_k\right)\\
&\geq \int_{D_2}|\nabla (\xi v_k)|^2-C\int_{D_2}\xi v_k^2-\alpha_* \delta_1\int_{D_2}\xi^2 v_k^2\\
&\quad\quad-\alpha_* (4\delta_1)^{-1}\int_{D_2}\xi^2+\epsilon_0e^{-\frac{2}{n}A}\frac{8}{n^3}\int_{D_2}\xi^2v_k^4\\
&\geq \int_{D_2}|\nabla (\xi v_k)|^2+\left(-C\delta_2+\epsilon_0e^{-\frac{2}{n}A}\frac{8}{n^3}\right)\int_{D_2}\xi^2v_k^4\\
&\quad\quad-\alpha_* \delta_1\int_{D_2}\xi^2 v_k^2-\left(C(4\delta_2)^{-1}+\alpha_* (4\delta_1)^{-1}\right)\text{vol}(D_2)\\
&\geq \int_{D_2}|\nabla (\xi v_k)|^2+\left(-C\delta_2+\epsilon_0e^{-\frac{2}{n}A}\frac{8}{n^3}-\alpha_* \delta_1\right)\int_{D_2}\xi^2v_k^2\\
&\quad\quad-\left(-C\delta_2+\epsilon_0e^{-\frac{2}{n}A}\frac{8}{n^3}\right)\int_{D_2}\xi^2-\left(C(4\delta_2)^{-1}+\alpha_* (4\delta_1)^{-1}\right)\text{vol}(D_2)
\end{aligned}
\end{equation}
Pick $\delta_1=\frac{2}{\alpha_*n^3}\epsilon_0e^{-\frac{2}{n}A}$ and $\delta_2=\frac{2}{Cn^3}\epsilon_0e^{-\frac{2}{n}A}$, we have
\begin{equation}
\begin{aligned}
 \int_{D_2}&|\nabla (\xi v_k)|^2+\frac{4}{n^3}\epsilon_0e^{-\frac{2}{n}A}\int_{D_2}\xi^2v_k^2\\
 &\leq \left(\frac{6}{n^3}\epsilon_0e^{-\frac{2}{n}A}+(\alpha_*^2+C^2)\epsilon_0^{-1}e^{\frac{2}{n}A}\frac{n^3}{8}\right)\text{Vol}(D_2)
\end{aligned}
\end{equation}
\qed

\begin{lemma}\label{lem3.3}
$u_k$ is locally uniformly $C^0$-bounded in the open set $M_-$ for $n\geq 2$ as $k\rightarrow +\infty$.
\end{lemma}
\proof
We choose cut-off functions $\eta\in C^\infty_0(M)$ at each step:
\begin{align}\label{3.11}
\eta(x)=
\begin{cases}
1& x\in B_{r_{i+1}}(x_0)\\
0& x\in M\setminus B_{r_ i}(x_0)
\end{cases}
,\ 0\leq\eta\leq1,\ |\nabla \eta|\leq \frac{C}{r_i-r_{i+1}},
\end{align}
where $B_{r}(x_0)=\{x\in M\ :\ dist(x_0, x)<r\}$, $r_i=\theta+\frac{\tau-\theta}{2^i}$, $\theta<r_i\leq\tau$ for $i=0,1,2, \cdots$, and $B_\tau(x_0)\subset M_-$.
Taking a test function $\varphi=\eta^2v_k^{a+1}\ (a\geq 0)$ in \eqref{3.5} yields
\begin{equation}\label{3.12}
\begin{aligned}
0&=\int_M\left(\nabla v_k\cdot \nabla (\eta^2v_k^{a+1})+\alpha_k\eta^2v_k^{a+1}-Se^{-\frac{2}{n}A}e^{\frac{2}{n}v_k}\eta^2v_k^{a+1}\right)d\mu_\omega\\
&=\int_M\left[\left(\frac{a}{2}+1\right)^{-1}(|\nabla (\eta v_k^{\frac{a}{2}+1})|^2-|\nabla \eta|^2v_k^{a+2})+\frac{a}{2}\eta^2v_k^a|\nabla v_k|^2\right]\\
&\quad\quad +\int_{B_{r_i}}\alpha_k\eta^2v_k^{a+1}-\int_{B_{r_i}}Se^{-\frac{2}{n}A}e^{\frac{2}{n}v_k}\eta^2v_k^{a+1}\\
&\geq \int_M\left[\left(\frac{a}{2}+1\right)^{-1}(|\nabla (\eta v_k^{\frac{a}{2}+1})|^2-|\nabla \eta|^2v_k^{a+2})\right]+\int_{B_{r_i}}\alpha_k\eta^2v_k^{a+1}\\
&\quad\quad-\int_{B_{r_i}}Se^{-\frac{2}{n}A}e^{\frac{2}{n}v_k}\eta^2v_k^{a+1}.
\end{aligned}
\end{equation}
By \eqref{3.11}, $B_{r_i}\subset M_-$ (i.e., $S|_{B_{r_i}}\leq -\epsilon_0<0$) and the fact $xe^x\geq -1$ for any $x\in\mathbb{R}$, we obtain
\begin{align}
&\int_M|\nabla \eta|^2v_k^{a+2}\leq \frac{C^2}{(r_i-r_{i+1})^2}\int_{B_{r_i}}v_k^{a+2},\\
&-\int_{B_{r_i}}\alpha_k\eta^2v_k^{a+1}\leq |\alpha_*|\int_{B_{r_i}}v_k^{a+1},\\
&\int_{B_{r_i}}Se^{-\frac{2}{n}A}e^{\frac{2}{n}v_k}\eta^2v_k^{a+1}\leq \frac{n}{2}\max_M|S| \int_{B_{r_i}} v_k^a.\label{3.15}
\end{align}
From \eqref{3.12}-\eqref{3.15}, it follows that
\begin{equation}
\begin{aligned}
\int_M|\nabla (\eta v_k^{\frac{a}{2}+1})|^2\leq &\frac{C^2}{(r_i-r_{i+1})^2}\int_{B_{r_i}}v_k^{a+2}+\left(\frac{a}{2}+1\right)|\alpha_*|\int_{B_{r_i}}v_k^{a+1}\\
&+\left(\frac{a}{2}+1\right)\frac{n}{2}\max_M|S| \int_{B_{r_i}} v_k^a.
\end{aligned}
\end{equation}
Using the H\"older inequality, we get
\begin{equation*}
\begin{aligned}
&\int_M v^a\leq \frac{a}{a+2}\int_M v^{a+2}+\frac{2}{a+2}\text{vol}(M),\\
&\int_M v^{a+1}\leq \frac{a+1}{a+2}\int_M v^{a+2}+\frac{1}{a+2}\text{vol}(M).
\end{aligned}
\end{equation*}
Thus,
\begin{equation}\label{3.17}
\int_M|\nabla (\eta v_k^{\frac{a}{2}+1})|^2\leq C\left(a+\frac{1}{2}+\frac{1}{(r_i-r_{i+1})^2}\right)\max\{\int_{B_{r_i}}v_k^{a+2},1\}.
\end{equation}
Applying Lemma \ref{lem2.2} to $f=\eta v_k^{\frac{a}{2}+1}$ gives
\begin{align}\label{3.18}
\left(\int_M (\eta v_k^{\frac{a}{2}+1})^{2\beta}\right)^{\frac{1}{\beta}}\leq C\left( \int_M|\nabla (\eta v_k^{\frac{a}{2}+1})|^2+\int_M(\eta v_k^{\frac{a}{2}+1})^2\right),
\end{align}
where $\beta=\frac{n}{n-1}$. By \eqref{3.17} and \eqref{3.18}, we obtain
\begin{align}\label{3.19}
\left(\int_{B_{r_{i+1}}}  v_k^{(a+2)\beta}\right)^{\frac{1}{\beta}}\leq C(a+2+\frac{1}{(r_i-r_{i+1})^2})\max\{\int_{B_{r_i}}v_k^{a+2},1\}.
\end{align}
Choose $a=a_i$ in \eqref{3.19} with $\beta^i=\frac{a_i}{2}+1$, we deduce that
\begin{equation}
\begin{aligned}
\max&\{\|v_k\|_{L^{2\beta^{i+1}}(B_{r_{i+1}})},1\}\\
&\leq C^{\frac{1}{2\beta^i}}\left( 2\beta^i+\frac{4^{i+1}}{(\tau-\theta)^2}\right)^{\frac{1}{2\beta^i}}\max\{\|v_k\|_{L^{2\beta^i}(B_{r_i})},1\}\\
&\leq C^{\frac{1}{2\beta^i}}4^{\frac{i+1}{2\beta^i}}\left(\frac{1}{2}+\frac{1}{(\tau-\theta)^2}\right)^{\frac{1}{2\beta^i}}\max\{\|v_k\|_{L^{2\beta^i}(B_{r_i})},1\}
\end{aligned}
\end{equation}
for any $i=0,1,2,\cdots$, since $\frac{\beta}{4}=\frac{n}{4(n-1)}<1$. By iteration, we get
\begin{equation}
\begin{aligned}
\max&\{\|v_k\|_{L^{2\beta^{i+1}}(B_{r_{i+1}})},1\}\\
&\leq 4^{\frac{1}{2}\sum_{k=0}^i\frac{k+1}{\beta^k}}\left(\frac{C}{2}+\frac{C}{(\tau-\theta)^2}\right)^{\frac{1}{2}\sum_{k=0}^i\frac{1}{\beta^k}}\max\{\|u_\lambda\|_{L^{2}(B_{\tau})},1\}. 
\end{aligned}
\end{equation}
Let $i\rightarrow +\infty$, then
\begin{align}
\max\{\|v_k\|_{C^{0}(B_{\theta})},1\}\leq C\max\{\|v_k\|_{L^{2}(B_{\tau})},1\}.
\end{align}
Using $v_k=u_k+A$ and Lemma \ref{lem3.2} yields
\begin{align}
\|u_k\|_{C^0(B_{\theta})}\leq C.
\end{align}
\qed

\begin{lemma}\label{lem3.4}
$e^{\frac{u_k}{n}}$ is uniformly $W^{1,2}(M)$-bounded for $n\geq 1$ as $k\rightarrow +\infty$.
\end{lemma}
\proof
Let $h\in C^\infty_0(M)$ be a cut-off function with $h<0$ on a open subset $D$ of $M_-$ and $h=0$ in $M\setminus D$. According to the fact (3) in Section 1, there exists a unique solution $f\in C^\infty(M)$ such that
\begin{align}
-\Delta^{Ch}_\omega f+\alpha_*=he^{\frac{2}{n}f}.
\end{align}
Set $w_k=u_k-f$, then
\begin{align}
-\Delta^{Ch}_\omega w_k+(\alpha_k-\alpha_*)=Se^{\frac{2}{n}u_k}-he^{\frac{2}{n}f}.
\end{align}
Multiplying the above equation by $e^{\frac{2}{n}w_k}$ and integrating on $M$ give
\begin{align}\label{3.27}
2n\int_M|\nabla e^{\frac{w_k}{n}}|^2+(\alpha_k-\alpha_*)\int_Me^{\frac{2}{n}w_k}= \int_M Se^{\frac{2}{n}(u_k+w_k)}-\int_M he^{\frac{2}{n}u_k}.
\end{align}
Pick $\varphi=e^{\frac{w_k}{n}}$ in \eqref{3.4}, we get
\begin{align}\label{3.28}
2\int_M|\nabla e^{\frac{w_k}{n}}|^2-\frac{4}{n}\int_M Se^{\frac{2}{n}(u_k+w_k)} \geq 0.
\end{align}
Substituting \eqref{3.28} into \eqref{3.27} and using $\alpha_k>\alpha_*$ yield
\begin{align}
\frac{3n}{2}\int_M|\nabla e^{\frac{w_k}{n}}|^2\leq -\int_D he^{\frac{2}{n}u_k}.
\end{align}
For $n=1$, by Lemma \ref{lemma2.4} and Lemma \ref{lem3.2}, we have 
\begin{align}\label{3.30}
\int_De^{2u_k}\leq C.
\end{align}
For $n\geq 2$, using Lemma \ref{lem3.3}, we have
\begin{align}
-\int_D he^{\frac{2}{n}u_k}\leq -\int_D he^{\frac{2}{n}\|u_k\|_{C^0(D)}}\leq C.
\end{align}
Hence, 
\begin{align}\label{3.32}
\int_M|\nabla e^{\frac{w_k}{n}}|^2\leq C
\end{align}
for any $n\geq 1$. Now we claim that $\|e^{\frac{w_k}{n}}\|_{L^2(M)}\leq C$. If it is not true, then there is a sequence $k_i\rightarrow +\infty$ such that $\lim_{i\rightarrow +\infty}\|e^{\frac{w_{k_i}}{n}}\|_{L^2(M)}=+\infty$. Let
\begin{align*}
g_i=\frac{e^{\frac{w_{k_i}}{n}}}{\|e^{\frac{w_{k_i}}{n}}\|_{L^2(M)}},
\end{align*}
then $\|g_i\|_{L^2(M)}=1$. Furthermore, it follows from \eqref{3.32} that $\lim_{i\rightarrow +\infty}\|\nabla g_i\|_{L^2(M)}=0$. Therefore, by passing to a subsequence, $g_i\rightharpoonup g_\infty$ weakly in $W^{1,2}(M)$, where $g_\infty\equiv C^*$ is a constant with $\|g_\infty\|_{L^2(M)}=1$. Here $\|g_\infty\|_{L^2(M)}=1$ is due to the Sobolev’s compact embedding theorem $W^{1,2}(M)\subset\subset L^2(M)$ and $\|g_i\|_{W^{1,2}(M)}\leq C$. For $n=1$, using \eqref{3.30}, we obtain that $\lim_{i\rightarrow +\infty}\|g_i\|_{L^2(D)}=0$. For $n\geq 2$, according to Lemma \ref{lem3.3}, we have $\|w_{k_i}\|_{L^\infty(D)}\leq C$, and thus $\lim_{i\rightarrow +\infty}\|g_i\|_{L^2(D)}=0$.  Hence, $C^*=0$ for any $n\geq 1$, which is a contradiction with $\|g_\infty\|_{L^2(M)}=1$. By $e^{\frac{1}{n}u_k}=e^{\frac{1}{n}w_k}e^{\frac{1}{n}f}$, we conclude that $e^{\frac{u_k}{n}}$ is uniformly $W^{1,2}(M)$-bounded for $n\geq 1$ as $k\rightarrow +\infty$.
\qed

\begin{lemma}\label{lem3.5}
$u_k$ is uniformly $W^{1,2}(M)$-bounded for $n\geq 1$ as $k\rightarrow +\infty$.
\end{lemma}
\proof
In terms of Lemma \ref{lem3.4}, we obtain that
\begin{equation}\label{3.32}
\begin{aligned}
\int_M |u_k|^p&=\int_{\{x\in M: -A\leq u_k(x)\leq 0\}} |u_k|^p+\int_{\{x\in M:  u_k(x)> 0\}} |u_k|^p\\
&\leq A^p\text{vol}(M)+\left(\frac{pn}{2}\right)^p\int_M e^{\frac{2}{n}u_k}\\
&\leq C
\end{aligned}
\end{equation}
for any $p>0$, where we have used the inequality $e^{px}>x^p\ (\forall x>0)$ and $C$ is a constant independent of $k$. Choose $\varphi=e^{\frac{u_k}{n}}$ in \eqref{3.4}, we get
\begin{align}\label{3.33}
\frac{1}{2n}\int_M|\nabla u_k|^2e^{\frac{2}{n}u_k}d\mu_\omega\geq\int_M Se^{\frac{4}{n}u_k} d\mu_\omega.
\end{align}
Multiplying \eqref{3.1} by $e^{\frac{u_k}{n}}$ and integrating on $M$ yield
\begin{align}\label{3.34}
\frac{2}{n}\int_M|\nabla u_k|^2e^{\frac{2}{n}u_k}+\alpha_k\int_Me^{\frac{2}{n}u_k}=\int_M Se^{\frac{4}{n}u_k}.
\end{align}
Combining \eqref{3.33} and \eqref{3.34}, we deduce that
\begin{align}
\frac{3}{2n}e^{-\frac{2}{n}A}\int_M|\nabla u_k|^2\leq \frac{3}{2n}\int_M|\nabla u_k|^2e^{\frac{2}{n}u_k}\leq -\alpha_k\int_Me^{\frac{2}{n}u_k}\leq -\alpha_*C,
\end{align}
where we have used $u_k>-A$, $\alpha\in (\alpha_*, \alpha^*)$ and Lemma \ref{lem3.4}. Therefore, 
\begin{align}
\int_M|\nabla u_k|^2\leq-\frac{2n}{3}e^{\frac{2}{n}A}\alpha_*C.
\end{align}
\qed

\begin{lemma}\label{lem3.7}
When $\dim_{\mathbb{C}}M=n\leq 2$, $u_k$ is uniformly $C^{2,\alpha}(M)$-bounded as $k\rightarrow +\infty$.
\end{lemma}
\proof
For $n=1$, using Lemma \ref{lem3.5} and Lemma \ref{lemma2.4}, we have
\begin{equation}\label{3.37}
\begin{aligned}
\int_Me^{pu_k}&= e^{p\bar{u}_k}\int_M e^{p(u_k-\bar{u}_k)}\\
&\leq e^{p\bar{u}_k}e^{\frac{p^2}{4\delta_1}\|\nabla u_k\|_{L^2(M)}^2}\int_M e^{\delta_1\left(\frac{u_k-\bar{u}_k}{\|\nabla u_k\|_{L^2(M)}}\right)^2}\\
&\leq C
\end{aligned}
\end{equation}
for any $p>0$, where $C$ is a constant independent of $k$. According to \eqref{3.32}, \eqref{3.37} and the elliptic $L^p$-estimate, we deduced that $u_k$ is uniformly $W^{2,p}(M)$-bounded for any $p>1$ as $k\rightarrow +\infty$. By Sobolev embedding theorem, $u_k$ is uniformly $C^{1,\alpha}(M)$-bounded for any $\alpha\in (0,1)$ as $k\rightarrow +\infty$. Using Schauder estimate, we obtain $u_k$ is uniformly $C^{2,\alpha}(M)$-bounded as $k\rightarrow +\infty$. 

For $n=2$, using the Sobolev embedding theorem $W^{1,2}(M)\subset L^4(M)\ (\dim_{\mathbb{R}}M=4)$, it follows from Lemma \ref{lem3.4} that 
\begin{align}\label{3.38}
\int_Me^{2u_k}\leq C.
\end{align}
By the elliptic $L^2$-estimate for \eqref{3.1} with $\alpha_*<\alpha_k<\alpha^*<0$, we obtain
\begin{equation}
\begin{aligned}
\|u_k\|_{W^{2,2}(M)}\leq C(\|u_k\|_{L^2(M)}+\max_M|S|\cdot\|e^{u_k}\|_{L^2(M)}).
\end{aligned}
\end{equation}
Using Lemma \ref{lem3.5} and \eqref{3.38} yields $\|u_k\|_{W^{2,2}(M)}$ is uniformly bounded. Then in terms of Lemma \ref{lemma2.4}, we have
\begin{equation}
\begin{aligned}
\int_Me^{pu_k}&= e^{p\bar{u}_k}\int_M e^{p(u_k-\bar{u}_k)}\\
&\leq e^{p\bar{u}_k}e^{\frac{p^2}{4\delta_1}\|\nabla^2 u_k\|_{L^2(M)}^2}\int_M e^{\delta_1\left(\frac{u_k-\bar{u}_k}{\|\nabla^2 u_k\|_{L^2(M)}}\right)^2}\\
&\leq C
\end{aligned}
\end{equation}
for any $p>0$, where $C$ is a constant independent of $k$. Following the same argument as the case $n=1$, it is easy to see that $u_k$ is uniformly $C^{2,\alpha}(M)$-bounded as $k\rightarrow +\infty$.
\qed

\section{An alternative proof of Ding-Liu's theorem}
In this section, inspired by an idea of \cite{[CL]}, we will give an alternative proof of Ding-Liu's theorem on compact Riemannian surfaces by using the $\sup +\inf$ inequality established by Brezis-Li-Shafrir \cite{[BLS]}. 
\begin{lemma}[\cite{[BLS]}]\label{lem4.1}
Let $V$ be a Lipschitz function satisfying 
\begin{align}
0<a\leq V(x)\leq  b<\infty
\end{align}
and $K$ a compact subset of a domain $\Omega\subset \mathbb{R}^2$. Then any solution $u$ of 
\begin{align}
-\Delta u=V(x)e^u, \quad x\in \Omega
\end{align}
satisfies 
\begin{align}
\sup_K u+\inf_K u\leq C(a,b,\|\nabla V\|_{L^\infty}, K,\Omega).
\end{align}
\end{lemma}

We recall the Ding-Liu's theorem \cite{[DL]} on prescribing sign-changing Gaussian curvatures:
\begin{theorem}\label{thm4.2}
Let $(M^n,\omega)$ be a compact Riemannian surfaces with Euler characteristic $\chi (M)<0$. Let $g_0$ be a nonconstant smooth function on $M$ with $\max_M g_0=0$. Then there exists a constant $\lambda^*\in (0,-\min_M g_0)$ such that 
\begin{enumerate}
\item If $\lambda\in (0, \lambda^*)$, then $g_0+\lambda\in PC(\omega)$;\\
 If $\lambda\in (\lambda^*,+\infty)$, then $g_0+\lambda\not\in PC(\omega)$.
 \item If $\lambda=\lambda^*$,  then $g_0+\lambda^*\in PC(\omega)$.
\end{enumerate}
Here $PC(\omega)$ denotes the set of $C^\infty(M)$ functions which are Gaussian curvatures of all $\tilde\omega\in \{\omega\}$.
\end{theorem}
\proof 
The existence of $\lambda^*$ and the assertion (1) is due to the implicit function theorem. Here we only prove the assertion (ii) by using the $\sup +\inf$ inequality. In Ding-Liu's proof, they have showed that for any $\lambda\in (0, \lambda^*)$, there exists a solution $u_\lambda\in C^\infty(M)$ of
\begin{align}\label{4.1}
-\Delta_\omega u_\lambda+s_0=(g_0+\lambda)e^{2u_\lambda}
\end{align}
satisfying 
\begin{align}
I_\lambda(u_\lambda)=\inf_{u\in X}I_\lambda(u)=\inf_{u\in X}\int_M \left(|\nabla u|^2+2s_0u-(g_0+\lambda)e^{2u}\right)d\mu_\omega
\end{align}
and 
\begin{align}\label{4.3}
0\leq\frac{d^2}{dt^2}|_{t=0}I_\lambda(u_\lambda+t\varphi)=2\int_M|\nabla \varphi|^2d\mu_\omega-4\int_M (g_0+\lambda)e^{2u_\lambda}\varphi^2 d\mu_\omega
\end{align}
for any $\varphi\in W^{1,2}(M)$, where $X=\{u\in W^{1,2}(M): u_1\leq u\leq u_2\ \text{a.e.\ on}\ M\}$, $u_1$ (Resp. $u_2$) is a smooth sub-solution (Resp. super-solution) of \eqref{4.1} with $u_1<u_2$. Here both $u_1$ and $u_2$ are related to $\lambda$. Moreover, they proved that there exists a constant $C>0$ independent of $\lambda$ such that $u_\lambda>-C$. 

Since $g_0+\lambda$ is sign-changing, then we can choose a small open subset $D_1\subset  M^+=\{x\in M: g_0(x)+\lambda>\epsilon_0\}$ for some $\epsilon_0>0$. Let $v$ be a solution of 
\begin{align}
\begin{cases}
-\Delta_\omega v=s_0\quad \text{in}\ D_1,\\
v=1\quad \text{on}\ \partial D_1.
\end{cases}
\end{align}
Set $w_\lambda=u_\lambda+v$, then
\begin{align}\label{4.5}
-\Delta_\omega w_\lambda=(g_0+\lambda)e^{-2v}e^{2w_\lambda},\quad \text{in}\ D_1.
\end{align}
Since $u_\lambda>-C$ on $M$ and $v$ is bounded in $D_1$, $w_\lambda$ is uniformly bounded below as $\lambda\rightarrow \lambda^*$. Since the metric is pointwise conformal to the Euclidean metric, applying Lemma \ref{lem4.1} to \eqref{4.5}, we obtain
\begin{align}
\sup_K w_\lambda\leq C(\epsilon_0, \lambda^*, \|\nabla g_0\|_{L^\infty}, K, D_1)-\inf_K w_\lambda<C
\end{align}
for any compact subset $K$ of $D_1$. Consequently, $w_\lambda$ is uniformly bounded in $K$ as $\lambda\rightarrow \lambda^*$, and thus $u_\lambda$ is also uniformly bounded in $K$.

Let $h\in C_0^\infty(M)$ be a nonpositive cut-off function such that $h<0$ in some open subset $D_2\subset M^+$ and $h\equiv 0$ in $M\setminus D_2$. Since $h\leq (\not\equiv) 0$, there exists a unique solution $\varphi\in C^\infty(M)$ of $-\Delta_\omega \varphi+s_0=he^{2\varphi}$. Set $f_\lambda=u_\lambda-\varphi$, then it satisfies 
$-\Delta_\omega f_\lambda=(g_0+\lambda)e^{2u_\lambda}-he^{2\varphi}$. Multiplying this equation by $e^{2f_\lambda}$ and integrating on $M$ give 
\begin{align}
2\int_M |\nabla e^{f_\lambda}|^2-\int_M (g_0+\lambda)e^{2(u_\lambda+f_\lambda)}=-\int_M h e^{2u_\lambda}.
\end{align}
Substituting \eqref{4.3} with $\varphi=e^{f_\lambda}$ into above equation yields
\begin{align}
\frac{3}{2}\int_M  |\nabla e^{f_\lambda}|^2\leq -\int_{D_2}h e^{2u_\lambda}\leq -(\inf_{D_2} h) e^{2\|u_\lambda\|_{L^\infty(D_2)}}\leq C,
\end{align}
where the last inequality holds because $u_\lambda$ is uniformly bounded in $D_2\subset M^+$. Following the proof of Lemma \ref{lem3.4}, it is easy to prove that $\|e^{f_\lambda}\|_{W^{1,2}(M)}\leq C$. By $e^{u_\lambda}=e^{f_\lambda}e^{\varphi}$, we obtain that $\|e^{u_\lambda}\|_{W^{1,2}(M)}\leq C$. According to the same argument of Lemma \ref{lem3.5} and Lemma \ref{lem3.7} with $n=1$, it is easy to see that $u_\lambda$ is uniformly $C^{2,\alpha}(M)$-bounded as $\lambda\rightarrow \lambda^*$. Therefore, there exists a subsequence $\{u_\lambda\}$ that converges to a function $u_{\lambda^*}$ in $C^{2,\alpha}(M)$, which is the desired solution of \eqref{4.1} for $\lambda=\lambda^*$.
\qed

\begin{remark}
For the case $\lambda=0$, following a similar argument above, we can also prove $g_0\in PC(\omega)$ under the same assumption as in Theorem \ref{thm4.2}. Indeed, let $\lambda\rightarrow 0^+$ in the above proof instead of $\lambda\rightarrow \lambda^*$, the conclusion can be obtained immediately.
\end{remark}

\section*{Appendix}
In this appendix, we will give a simple proof of the key Lemma \ref{lem3.3} in the proof of Theorem \ref{thm1.1} by using the maximum principle. The advantage of this approach is that we don't need to assume that the manifold $(M^n, \omega)$ is balanced.
\begin{proposition}
Let $(M^n, \omega)$ be a compact Hermitian manifold with $\dim_{\mathbb{C}}M=n\geq 1$. Suppose that $u_k\in C^\infty(M)$ satisfy
\begin{align}\label{4.12}
-\Delta^{Ch}_\omega u_k+\alpha_k=Se^{\frac{2}{n}u_k}
\end{align}
where $\{\alpha_k\}$ is a sequence of numbers with $\alpha_*<\alpha_k<\alpha^*<0$ and $\lim_{k\rightarrow +\infty}\alpha_k=\alpha_*$, $\alpha^*<0$ is a fixed constant, $S$ is a sign-changing smooth function on $M$ with $\int_M S d\mu_\omega<0$. For any compact subset $K$ of $M_-=\{x\in M: S(x)<0\}$, there exists a uniform constant $C_1=C_1(S, \alpha_*, K)>0$ such that
\begin{align}
\sup_Ku_k <C_1.
\end{align}
Combining with Lemma \ref{lemma 3.3}, we get 
\begin{align}
\|u_k\|_{C^0(K)}\leq C_2
\end{align}
for some positive constant $C_2$ independent of $k$.
\end{proposition}
\proof Set 
\begin{align}
v_k=e^{\frac{2}{n}u_k}.
\end{align}
By \eqref{4.12}, we obtain
\begin{align}\label{4.15}
-\Delta^{Ch}_\omega v_k+\frac{|\nabla v_k|^2}{v_k}+\frac{2}{n}\alpha_kv_k=\frac{2}{n}Sv_k^2.
\end{align}
We choose a cutoff function $\varphi\in C^\infty_0(M_-)$ with $0\leq \varphi\leq 1$ on $M$ and $\varphi\equiv 1$ on $K$. Set
\begin{align}
\tilde{v}_k=\varphi^2e^{\frac{2}{n}u_k},
\end{align}
then there exists a point $x_0\in M_-$ such that $\tilde{v}_k(x_0)=\max_M \tilde{v}_k>0$. By the maximum principle, we have
\begin{align}\label{4.17}
0=\nabla \tilde{v}_k(x_0)=\nabla(\varphi^2v_k)=2\varphi v_k \nabla \varphi +\varphi^2\nabla v_k,
\end{align}
and
\begin{equation}\label{4.18}
\begin{aligned}
0&\geq \Delta^{Ch}_\omega \tilde{v}_k(x_0)\\
&=\Delta^{Ch}_\omega (\varphi^2v_k)\\
&=2|\nabla \varphi|^2v_k+2\varphi v_k\Delta^{Ch}_\omega\varphi +4\varphi\nabla\varphi \cdot \nabla v_k+\varphi^2\Delta^{Ch}_\omega v_k.
\end{aligned}
\end{equation}
Substituting \eqref{4.15} and \eqref{4.17} into \eqref{4.18} yields at $x_0$
\begin{equation}
\begin{aligned}
0&\geq2|\nabla \varphi|^2v_k+2\varphi v_k\Delta^{Ch}_\omega\varphi +4\varphi\nabla\varphi \cdot \nabla v_k+\varphi^2\Delta^{Ch}_\omega v_k\\
&=2|\nabla \varphi|^2v_k+2\varphi v_k\Delta^{Ch}_\omega\varphi +4\varphi\nabla\varphi \cdot (-2\varphi^{-1}v_k \nabla \varphi)\\
&\quad\quad+\varphi^2\left( \frac{|-2\varphi^{-1}v_k \nabla \varphi|^2}{v_k}+\frac{2}{n}\alpha_kv_k-\frac{2}{n}Sv_k^2\right)\\
&=-2|\nabla \varphi|^2v_k+2\varphi v_k\Delta^{Ch}_\omega\varphi +\varphi^2\frac{2}{n}\alpha_kv_k-\frac{2}{n}S\varphi^2v_k^2.
\end{aligned}
\end{equation}
Thus, there is a positive constant $C=C(K, \|\varphi\|_{C^2(M)}, \alpha_*)$ such that
\begin{equation}
\begin{aligned}
-\frac{2}{n}\max_{\text{supp}\varphi}S\cdot\max_M \tilde{v}_k&\leq-\frac{2}{n}S(x_0)\max_M \tilde{v}_k\\
&=-\frac{2}{n}S(x_0)(\varphi^2v_k)(x_0)\\
&\leq 2|\nabla \varphi|^2-2\varphi \Delta^{Ch}_\omega\varphi -\varphi^2\frac{2}{n}\alpha_*\\
&\leq C,
\end{aligned}
\end{equation}
which implies that for any $x\in K$
\begin{align}
e^{\frac{2}{n}u_k(x)}= (\varphi^2v_k)(x)=\tilde{v}_k(x)\leq \max_M \tilde{v}_k\leq -\frac{n}{2}(\max_{\text{supp}\varphi}S)^{-1}C,
\end{align}
where we have used $\varphi|_K\equiv 1$ and $K\subset M_-$. Therefore, we get $\sup_Ku_k <C$ for some constant $C>0$.
\qed

\section*{Statements and Declarations}
\textbf{Data availability.} No data was used for the research described in the article.

\textbf{Conflict of interest.}  The author states that there is no conflict of interest. The author has no relevant financial or non-financial interests to disclose.

\bigskip

Weike Yu

School of Mathematical Sciences, 

Ministry of Education Key Laboratory of NSLSCS,

Nanjing Normal University,

Nanjing, 210023, Jiangsu, P. R. China,

wkyu2018@outlook.com

\bigskip


\begin{thebibliography}{99}
\bibitem{[ACS]} D. Angella, S. Calamai, C. Spotti, On the Chern-Yamabe problem, Math. Res. Lett., 24, 645-677 (2017)

\bibitem{[BLS]} H. Brezis, Y. Li, I. Shafrir, A $\sup+\inf$ inequality for some nonlinear elliptic equation involving exponential nonlinearity, J. Funct. Anal., 115, 344-358 (1993).

\bibitem{[CZ]} S. Calamai, F. Zou, A note on Chern-Yamabe problem, Differential Geom. Appl., 69, 101612 (2020).

\bibitem{[CL]} W. Chen, C. Li, Gaussian curvature in the negative case, Proc. Amer. Math. Soc., 131(3), 741-744 (2002).

\bibitem{[DL]} W.-Y. Ding, J.-Q. Liu, A note on the problem of prescribing Gaussian curvature on surfaces, Trans. Amer. Math. Soc., 347(3), 1059-1066 (1995) 

\bibitem{[Doa]} A. Doan, Adiabatic limits and Kazdan–Warner equations. Calc. Var., 57, 124 (2018).
\bibitem{[Fus]} E. Fusi, The prescribed Chern scalar curvature problem, J. Geom. Anal., 32, 187 (2022).

\bibitem{[Fon]} L. Fontana, Sharp borderline Sobolev inequalities on compact Riemannian manifolds, Comm. Math. Helv., 68: 415-454 (1993).

\bibitem{[Gau1]} P. Gauduchon, Le th\'eor\`eme de l'excentricit\'e nulle, C. R. Acad. Sci. Paris S\'er., A-B 285(5), A387–A390 (1977).

\bibitem{[Gau]} P. Gauduchon, La 1-forme de torsion d'une vari\'et\'e hermitienne compacte, Math. Ann., 267(4), 495–518 (1984).

\bibitem{[Ho]}P.-T. Ho, Results related to the Chern-Yamabe flow, J. Geom. Anal., 31, 187-220 (2021).

\bibitem{[KW]} J. L. Kazdan, F. W. Warner, Curvature functions for compact 2-manifolds, Ann. of Math., 99(1): 14-47 (1974). 

\bibitem{[LM]} M. Lejmi, A. Maalaoui, On the Chern-Yamabe flow, J. Geom. Anal., 28, 2692-2706 (2018).

\bibitem{[Yu1]} W. Yu, Prescribed Chern scalar curvatures on compact Hermitian manifolds with negative Gauduchon degree, J. Funct. Anal., 285(2), 109948 (2023).

\bibitem{[Yu2]} W. Yu, A note on Kazdan-Warner type equations on compact Riemannian manifolds, Nonlinear Anal. TMA, 246, 113596 (2024).




\end{thebibliography}
\end{document}